\documentclass[11pt,a4wide]{article}

\usepackage{multirow}
\usepackage{algorithmic}
\usepackage{amsmath}
\usepackage{amsfonts}
\usepackage{amsthm}
\usepackage[english]{babel} 
\usepackage{dsfont}
\usepackage{amssymb}
\usepackage{graphicx}        
\usepackage{cite}                            
\usepackage{subcaption}
\usepackage{epstopdf}
\usepackage{epsfig}
\usepackage{subfloat}
\usepackage{float}
\usepackage[thinlines]{easytable}
\usepackage{array}
\usepackage{mathtools}
\usepackage{stmaryrd}
\usepackage{bm}
\usepackage{hyperref}
\hypersetup{
    colorlinks=true,
    linkcolor=blue,
    filecolor=magenta,      
    urlcolor=cyan,
    pdftitle={Overleaf Example},
    pdfpagemode=FullScreen,
    }

\usepackage{xcolor,colortbl}
\definecolor{g}{rgb}{0.68,1,0.18}

\usepackage{booktabs}
\usepackage{url}
\usepackage[font=small,labelfont=bf,textfont=it,indention=.1cm,width=1.1\textwidth]{caption}

\usepackage[capitalise]{cleveref}
\usepackage{comment}

\allowdisplaybreaks
\numberwithin{equation}{section}

\newtheorem{definition}{Definition}[section]

\usepackage[linesnumbered,algoruled,boxed,lined]{algorithm2e}

\def\be{\begin{equation}}
\def\ee{\end{equation}}
\def\bea{\begin{eqnarray}}
\def\eea{\end{eqnarray}}

\def\Y{F}

\title{Repulsion dynamics for uniform Pareto front approximation in multi-objective optimization problems}

\author{Giacomo Borghi\footnote{RWTH Aachen University, Institute for Geometry and Applied Mathematics, Aachen, Germany (borghi@eddy.rwth-aachen.de)}}

\begin{document}
\maketitle

\begin{abstract} 
Scalarization allows to solve a multi-objective optimization problem by solving many single-objective sub-problems, uniquely determined by some parameters. In this work, we propose several adaptive strategies to select such parameters in order to obtain a uniform approximation of the Pareto front. This is done by introducing a heuristic dynamics where the parameters interact through a binary repulsive potential. The approach aims to minimize the associated energy potential which is used to quantify the diversity of the computed solutions. A stochastic component is also added to overcome non-optimal energy configurations. Numerical experiments show the validity of the proposed approach for bi- and tri-objectives problems with different Pareto front geometries.
\end{abstract}

{\bf Keywords}:  multi-objective optimization, gradient-free methods, heuristic algorithms, 

potential based diversity measure, scalarization


\section{Introduction}
In this work, we are interested in multi-objective optimization problems of the form
\begin{equation}
\textup{minimize} \;f(x) := (f_1(x), \dots, f_m(x))^\top \quad \textup{subject to} \quad x \in \mathcal{D}
\label{eq:mop}
\end{equation}
where $f \in \mathcal{C}(\mathcal{D}, \mathbb{R}^m)$ is a vector of $m\geq2$ objectives and $\mathcal{D} \subset\mathbb{R}^{\textup{d}}$ is the admissible domain.
Problems of type \eqref{eq:mop} often arise in applications whenever two or more objectives need to be minimized (or maximized), but there is no point $x^* \in \mathcal{D}$ which minimizes all of them at the same time. Therefore, the notion of \textit{Pareto optimality} is introduced to define solutions to \eqref{eq:mop}. We recall it from \cite{jahn2004vector} for completeness.

\begin{definition} A point $x^*$ is \textit{weakly Pareto optimal} if there is no other $x \in \mathcal{D}$ such that $f_\ell(x) < f_\ell(x^*)$ for all $\ell = 1, \dots, m$. Similarly, a point $x^*$ is \textit{strongly Pareto optimal} if there is no other $x \in \mathcal{D}$ such that $f_\ell(x) \leq f_\ell(x^*)$ for all $\ell = 1, \dots, m$ and such that $f(x) \neq f(x^*)$.
\end{definition}

Let $\mathcal{X}\subset{\mathbb{R}^\textup{d}}$ be the set of strongly Pareto optimal points and $\mathcal{F}:=f(\mathcal{X})$ its image on the objectives space, usually called \textit{Pareto front} \cite{jahn2004vector}. 
Solving \eqref{eq:mop} requires to find $N \in \mathbb{N}$ strongly Pareto optimal solutions which should be as \textit{diverse} as possible from one another to better describe the (potentially uncountable) set of optimal solutions \cite{book2005mop}. In terms of Pareto front $\mathcal{F}$, this translates into finding a $N$ solutions $x^{*,i}$ such that $\{f(x^{*,i}) \}_{i=1}^N$ are uniformly distributed over $\mathcal{F}$. There is no uniquely recognized way to measure uniformity, but different measures have been used in the literature, such as \textit{hypervolume contribution} \cite{zitzler1998multi}, \textit{crowding distance} \cite{deb2002nsga2} and, recently, the \textit{Riesz s-energy} \cite{coello2021overview, vega2021towards}.

A common way to address \eqref{eq:mop} is trough scalarization, as it allows to find weakly Pareto optimal solutions by solving a set of parametrized single-objective sub-problems \cite{jahn2004vector}. A classical example is the weighted sum approach where all the objectives are sum together, with certain weights, and then the sum is minimized. Scalarization, though, only addresses the Pareto optimality of the solutions and not their uniformity over the front $\mathcal{F}$. Which sub-problems lead to a good approximation of $\mathcal{F}$ is clearly not known in advance and so several criteria to select the sub-problems have been proposed in the literature, see e.g. \cite{braun2015preference,coello2017hyper}.


With the aim of finding sub-problems leading to a uniform approximation of $\mathcal{F}$, we propose heuristic adaptive strategies which rely on energy-based uniformity measures. Such measures have gained popularity in the multi-objective evolutionary optimization community thanks to their scalability and flexibility \cite{coello2020survey}. The method introduces a short-range repulsion dynamics between the parameters associated to the sub-problems and it aims to find a low-energy configuration over the Pareto front $\mathcal{F}$.
Such an approach was first introduced in \cite{borghi2022adaptive} where it was use together with a Consensus-Based Optimization (CBO) method \cite{pinnau2017consensus} to simultaneously solve $N$ different sub-problems. Here, we place the adaptive strategies in a more general framework by considering the situation where an arbitrary auxiliary method for solving the sub-problems is available. We also improve the adaptation process by adding stochasticity to the dynamics and perform novel numerical experiments with tri-objectives test problems.

The framework we consider and the adaptive strategies are presented in Section \ref{s2}. In Section \ref{s3} we perform numerical experiments by coupling the dynamics with a CBO solver and show the validity of the proposed approach for $m=2$ and $m=3$. Some final remarks and possible future research directions follows in Section \ref{s4}.

%
%
%
%
%
%
%
%

\section{Repulsive dynamics}
\label{s2}

Scalarization functions are a fundamental tool in multi-objective optimization as they allow to break \eqref{eq:mop} into several single-objective optimization sub-problems. A popular choice consist of the family of weighted $\ell_p$ semi-norms \cite{book2005mop,jahn2004vector}, $p \in [1, \infty)$
\begin{equation}
\textup{minimize}\;\; S(f(x), w):= \left (  \sum_{\ell=1}^m\, w_\ell |f_\ell(x) - \hat y_\ell |^p \right)^{\frac1p} \quad \text{subject to}\quad x \in \mathcal{D}
\label{eq:sub}
\end{equation}
where $ \hat y $ is an \textit{ideal} point, that it is, a collection of lower bounds for the objectives: $\hat y_{\ell} \leq f_\ell(x)$ for all $x \in \mathcal{D}$ and $\ell = 1, \dots, m$. The scalar sub-problems \eqref{eq:sub} are parametrized by a vector of weights $w = (w_1, \dots, w_m)$ belonging to the $(m-1)$-unit, or probability, simplex
\begin{equation*} \Omega : = \Big\{w \in \mathbb{R}^m_{\geq 0}\; | \; \sum_{\ell = 1}^m \,w_\ell = 1 \Big \}\,.\end{equation*}
We note that for $p = 1$ one obtains the weighted sum approach, while by extending the definition to $p=\infty$ by continuity, one obtains the Chebyschev scalarization strategy \cite{book2005mop}. If $x^*$ is a solution to \eqref{eq:sub} for some $w\in \Omega$, then $x^*$ is a weakly Pareto optimal point, and the reverse statement is also true if $p=\infty$ \cite{jahn2004vector}.


In the following, $p\in [1, \infty]$ is arbitrary but fixed and we assume to have an auxiliary solver $\xi$ which, given a parameter $w\in \Omega$, returns an approximation of solution to the correspondent sub-problem:
\[ \xi(w) \approx x^* \in \underset{x \in \mathcal{D}}{\textup{argmin}}\; S(f(x),w) \,.  \] 
Let $\{W_0^i\}_{i=1}^N \subset \Omega$ be an initial selection of sub-problems where the parameters are uniformly distributed over the simplex $\Omega$ and let $\Y_0^i := f(\xi(W^i_0))$ be the image of the solutions computed by the solver. It is well-known that $\{ \Y_0^i\}_{i=1}^N$ may not be well-distributed over the front $\mathcal{F}$ \cite{Das1997}. Since energy-based measure have shown to be a valid tool to quantify the concept of uniformity \cite{coello2021overview}, this means that $\{\Y^i_0\}_{i=1}^N$ may correspond to a high energy configuration of points over the objectives space $\mathbb{R}^m$. To be more precise, let $\mu \in \mathcal{P}(\mathbb{R}^m)$ be collection of points given as a probability measure. We quantify the uniformity of such distribution by means of a potential energy
\begin{equation} 
\mathcal{U}[\mu] := \iint U(y - z)\, d\mu(z)\,d \mu(y)
\label{eq:pot}
\end{equation}
where $U\in \mathcal{C}^1(\mathbb{R}^m\setminus \{0\})$ describes a short-range repulsion, like the Riesz potential $U(z) = \|z\|^{-s}, s>0$ or the Morse potential $U(z) = e^{-C\|z\|}, C>0$. 
In these terms, $\{F_0^i\}_{i=1}^N$ are not well-distributed if $\mathcal{U}[\mu^N_0] \gg \inf \{ \mathcal{U}[\nu]\,|\, \nu \in \mathcal{P}(\mathcal{F})\}$ with $\mu^N_0$ being the corresponding empirical measure $\mu_0^N = 1/N\sum_{i} \delta_{Y_0^i}$.

To reach a lower energy configuration, we first propose an adaptive process which is based on the vector field associated to \eqref{eq:pot} in the case $m=2$. Let $k = 0,1,2,\dots$ indicate the iteration step and $\Y^i_k$ be always the solutions computed by the auxiliary solver, $\Y_k^i := f(\xi(W^i_k))$. For a fixed step length $\tau>0$ we have 
\begin{equation}
W^i_{k+1} = \Pi_{\Omega} \bigg( W^i_k  + \frac{\tau}N \sum_{j=1}^N \nabla U(\Y_k^i - \Y_k^j)  \bigg) \qquad \textup{for}\quad i=1,\dots,N
\label{eq:old}
\end{equation}
where the gradient is taken with respect to the first argument and $\nabla U(0):=0$. The operator $\Pi_{\Omega}$ performs the projection to the simplex and it is well-defined thanks to the convexity of $\Omega$. Dynamics \eqref{eq:old} was first proposed and studied in \cite{borghi2022adaptive} in the more specific context of CBO methods. As shown in \cite{borghi2022adaptive}, this adaptive strategy exploits the (implicit) parametrization of the front given by the Chebyschev scalarization approach to simulate a gradient flow dynamics in the image space. This can be done for $m=2$ thanks to the intuitive relation between $\Omega$ and $\mathcal{F}$ and, in particular, because $W^i_k$ can only move in one direction given by $(1,-1)^\top$. 

When $m>2$ such relation between  $\Omega$ and $\mathcal{F}$ becomes more involved. With the same spirit, though, we define a repulsive dynamics in $\Omega$ which takes into account the forces computed in the objectives space. Let $U$ be radially symmetric and given by $U(z) = \overline U(\|z\|)$ for some $\overline U \in \mathcal{C}^1(\mathbb{R}_{>0})$, the dynamics reads
\begin{equation}
W^i_{k+1} = \Pi_{\Omega} \bigg( W^i_k  - \frac{\tau}N \sum_{j=1}^N  D\overline U(\|\Y_k^i - \Y_k^j\|)\frac{W_k^i - W_k^j}{\|W_k^i - W_k^j \|}  \bigg) \qquad \textup{for}\quad i=1,\dots,N
\label{eq:new}
\end{equation}
with $D\overline U$ being the derivative of $\overline U$ and $D\overline U(0) = 0$.
Here, the displacement direction is simply determined by the difference between $W^i_k$ and $W_k^j$, but the magnitude and the sign are given by the forces computed in the objective space. Indeed, in these settings, we can also have a long-range attraction between the parameters. Dynamics \eqref{eq:new} was only suggested in \cite{borghi2022adaptive}: in the next section we will compare it with \eqref{eq:old} and also validate it for $m=3$.

The described adaptive processes can be seen as a system of particles which evolve according to a mean-field type binary interaction over a bounded domain \cite{fatecau2017swarm, carrillo2016nonlocal}. By studying the corresponding many-particle limit, the authors in \cite{fatecau2017swarm} show that the system may evolve to non-optimal equilibria where a certain amount of particles, or mass, is concentrated at the boundary. This is due to the repulsive behavior and it suggests that, with dynamics \eqref{eq:old} and \eqref{eq:new}, some parameters $W^i_k$ may reach the boundary of the simplex and remain \textit{stuck} there, leading to a non-optimal configuration in terms of potential energy. Inspired by \cite{fatecau2019diffusion} where a diffusive component is added to overcome this effect, we propose a third dynamics where we add Gaussian noise to \eqref{eq:new}:
\begin{equation}
\begin{cases}
W^i_{k+1/2} &= \Pi_{\Omega} \left( W^i_k  - \frac{\tau}N \sum_{j=1}^N  D\overline U(\|\Y_k^i - \Y_k^j\|)\frac{W_k^i - W_k^j}{\|W_k^i - W_k^j \|}  \right) \\
W^i_{k+1} &=  \Pi_{\Omega} \left(  W_{k+1/2}^i + \zeta B_{k}^i \right)
\end{cases}
 \qquad \textup{for}\quad i=1,\dots,N
\label{eq:newnoise}
\end{equation}
where $\zeta>0$ is a parameter controlling the variance and $B_k^i \in \mathbb{R}^m$ are vectors independently sampled from the standard normal distribution, $B_k^i \sim \mathcal{N}(0,I_m)$. We note that the random component is added after a first projection step and not together with the repulsive forces. This does not make a difference if $W^i_k$ belongs to the interior of the simplex, but it increases the chances of the random component to be effective if $W^i_k$ belongs to the simplex boundary - exactly where aggregates of particles are likely to appear.










%
%
%
%
%
%

\section{Numerical experiments}
\label{s3}

In the considered settings, we require and auxiliary solver $\xi$ which returns a solution to every scalar sub-problem of the form \eqref{eq:sub}. As already mentioned, $\xi$ can be any suitable single-objective optimization method, but one should aware of the computational cost of solving several sub-problems of the type \eqref{eq:sub}.
To reduce the computational effort, we suggest to couple the adaptive strategy with a population-based solver which is able to solve many sub-problems at the same time, like the MOEA/D algorithm \cite{zhang2008moead}. We refer to \cite{deb2001multi,reyes2006pso} for more details on population-based meta-heuristics algorithms. 

In our experiments, we couple the the adaptive strategy with the Multi-objective Consensus-Based Algorithm (M-CBO) \cite{borghi2022multi}, as done in \cite{borghi2022adaptive}. 
A precise presentation of the M-CBO iterative dynamics is beyond the purpose of this work. We will outline the main properties, while referring to \cite{borghi2022adaptive} (and the references therein) for more details on CBO methods.
In M-CBO, a population of agents explore the search space according to some deterministic and stochastic rules. 
A solution guess for every sub-problem $w \in \Omega$ is available and it is given by a weighted sum of the agents position. Therefore, we exploit this to update the parameters $\{W_k^i\}_{i=1}^N$ according to one of the proposed strategies. The parameters adaptation is performed only every $t_k=50$ iteration of the M-CBO solver to save computational cost. Algorithm \ref{alg} presents in the details the described procedure. We note that, differently to \cite{borghi2022adaptive}, Algorithm \ref{alg} uses $n = 20$ agents per sub-problem.

\begin{algorithm}[t]
\caption{} 
\label{alg}
\begin{algorithmic}[1]
\STATE{Set adaptive process parameters: $\tau, \zeta$, $t_k = 50$}
\STATE{Set CBO parameters: $\alpha = 10^{5}, \lambda =1, \sigma =1 , \Delta t = 10^{-2}$}
\STATE{Select the sub-problems $\{ W_0^i\}_{i=1}^N $ uniformly in $\Omega$. $N=15$ for $m=2$, $N=66$ for $m=3$.}
\STATE{Initialize the set: $X^{i,h}_0$ uniformly sampled from $\mathcal{D} = [0,1]^\textup{d}$, $i=1, \dots, N$, $h = 1, \dots, n=20$. }
\STATE{$k=0$}
\FOR{$s=0,\dots,s_{max} =200t_k$}
	\STATE{$\omega^{i,h} = \exp(-\alpha S(f(X_s^{i,h}),W_k^i))$}
	\STATE{$Y_s^{\alpha,i} = \sum_{i,h}X_s^{i,h}\omega^{i,h}/ \sum_{i,h}\omega^{i,h}$}
	\STATE{$X_{s+1/2}^{i,h} =  X_s^{i,h} + \lambda \Delta t(Y^{\alpha,i}_s -  X_s^{i,h})  + \sigma \sqrt{\Delta t}\|Y^{\alpha,i}_s -  X_s^{i,h}\|Z_{t}^{i,h}\quad $ with $\;  Z_{t}^{i,h} \sim \mathcal{N}(0,I_{\textup{d}})$}
	\STATE{$X_{s+1}^{i,h} = \Pi_{[0,1]^d}( X_{s+1/2}^{i,h})$}
       \IF{$\mod(t,t_k) = 0$ }
       \STATE{$\xi(W_k^i) = Y_s^{\alpha,i}$}
	\STATE{$F_k^i = f(\xi(W_k^i))$}
	\STATE{update $W_k^i$ according to \eqref{eq:old} or \eqref{eq:new}, \eqref{eq:newnoise}}
	\STATE{$k = k+1 $}
	\ENDIF
\ENDFOR
\RETURN $\{ Y^{\alpha,i}_s \}_{i=1}^N$  $\{ F^{i}_k \}_{i=1}^N$
\end{algorithmic}
\end{algorithm}

 \begin{figure}
\centering
 \includegraphics[trim = 2.6cm 5cm 2.1cm 0, clip, width=0.5\linewidth]{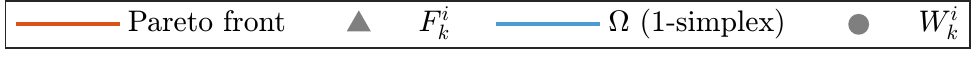}
\smallskip

\includegraphics[width=0.21\linewidth]{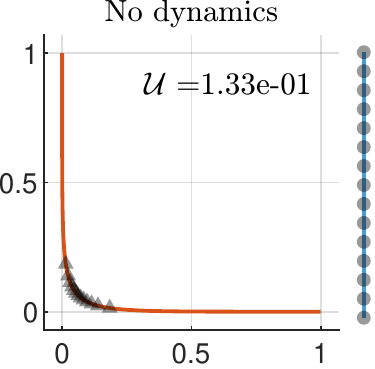}~a) 
\hfil
\includegraphics[width=0.21\linewidth]{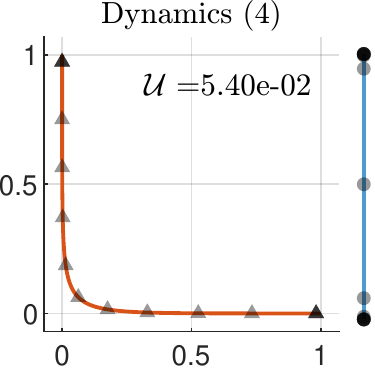}~b)
\hfil
\includegraphics[width=0.21\linewidth]{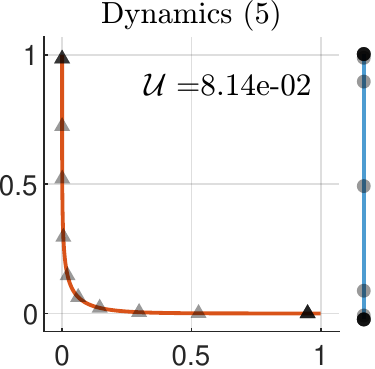}~c)
\includegraphics[width=0.21\linewidth]{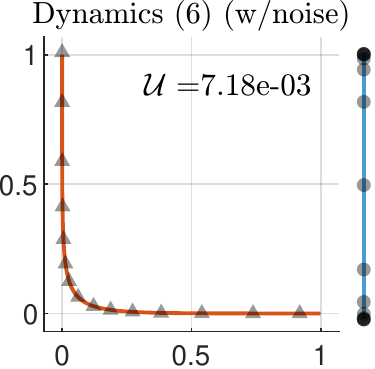}~d)
\caption{Bi-objective Lam\'e problem with $\gamma = 0.25.$}
\label{fig:1}
\end{figure}

We test Algorithm \ref{alg} with the Lam\'e problems \cite{emmerich2007lame}, where a parameter $\gamma>0$ controls the front curvature, and with the inverted DTLZ1 problem \cite{ishi2017} in order to cover a different variety of front shapes. As scalarization functions, we use the Chebyschev weighted semi-norms ($p = \infty$). We set $\textup{d}=m$, the strength of the interaction, if present, to $\tau = 10^{-2}$ and use the Morse potential $U(z) = e^{-C\|z\|}$ with $C=30$. The parameter controlling the magnitude of the random component is set to $\zeta=10^{-9}$ when $m=2$ and to $\zeta=10^{-6}$ when $m=3$. We update $W_k^i$ for 200 steps maximum. The remaining parameters controlling the CBO dynamics are specified in Algorithm \ref{alg}.

Fig.\ref{fig:1} shows the final distribution of $N=15$ optimal points computed by the algorithm with the different adaptive strategies proposed for a bi-objective problem. Consistently with \cite{borghi2022adaptive}, dynamics \eqref{eq:old} leads to well-distributed points over the front as well as when dynamics \eqref{eq:new} is used. We note 
how several points are concentrated at the extrema of the front: this is because the corresponding points over the simplex are stuck, due to the projection step, at the extrema of the simplex. As expected, adding a stochastic components mitigates this effect and the final configuration attains a lower potential value, see Fig.\ref{fig:1}d.

\begin{figure}
\centering
 \includegraphics[trim = 2.6cm 5cm 2.1cm 0, clip, width=0.5\linewidth]{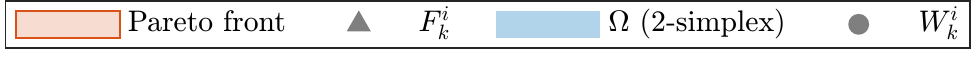}
\smallskip

\includegraphics[height=0.21\linewidth]{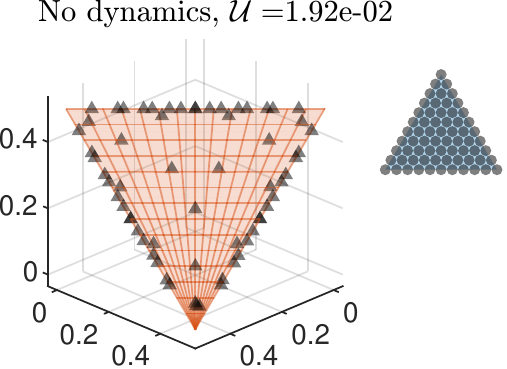}~a) 
\hfil
\includegraphics[height=0.21\linewidth]{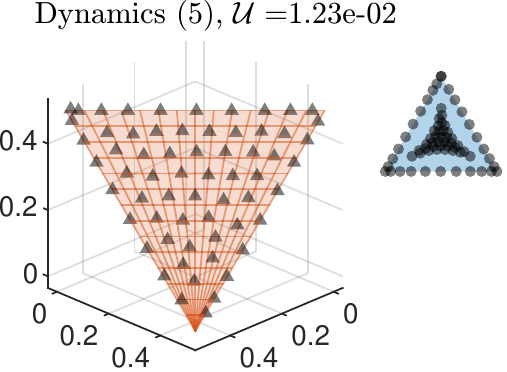}~b)
\hfil
\includegraphics[height=0.21\linewidth]{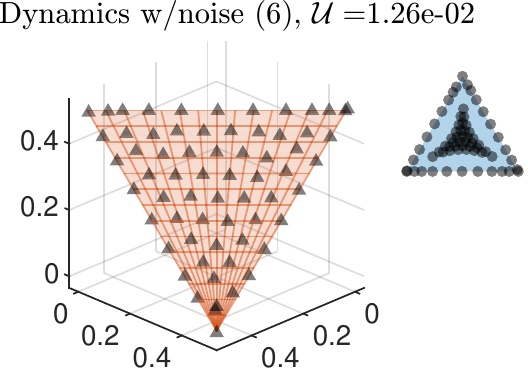}~c)
\caption{Tri-objective inverted DTLZ1 problem.}
\label{fig:2}
\end{figure}

\begin{figure}
\centering
\includegraphics[height=0.223\linewidth]{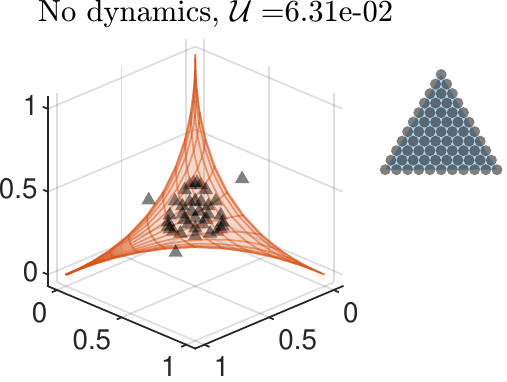}~a) 
\hfil
\includegraphics[height=0.223\linewidth]{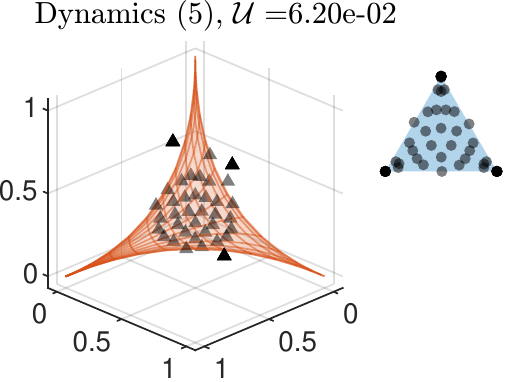}~b)
\hfil
\includegraphics[height=0.223\linewidth]{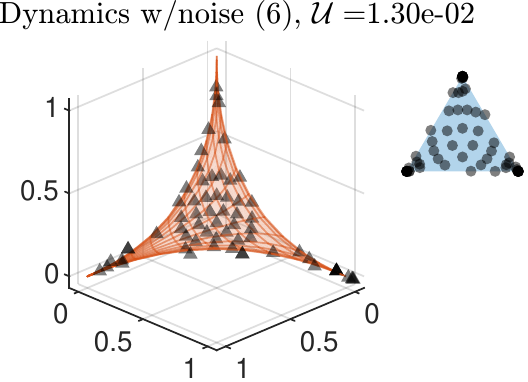}~c)
\caption{Tri-objective Lam\'e problem with $\gamma = 0.5.$}
\label{fig:3}
\end{figure}

\begin{figure}
\includegraphics[height=0.223\linewidth]{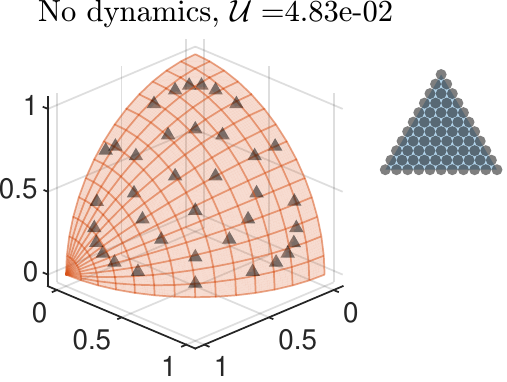}~a) 
\hfil
\includegraphics[height=0.223\linewidth]{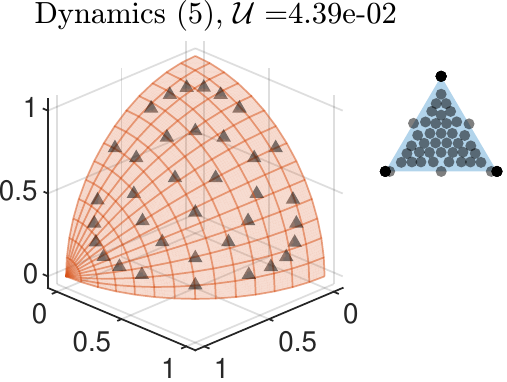}~b)
\hfil
\includegraphics[height=0.223\linewidth]{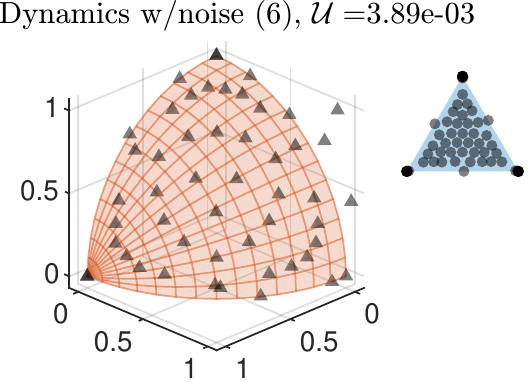}~c)
\caption{Tri-objective Lam\'e problem with $\gamma = 2.$}
\label{fig:4}
\end{figure}

\begin{figure}[t]
\includegraphics[width=\linewidth]{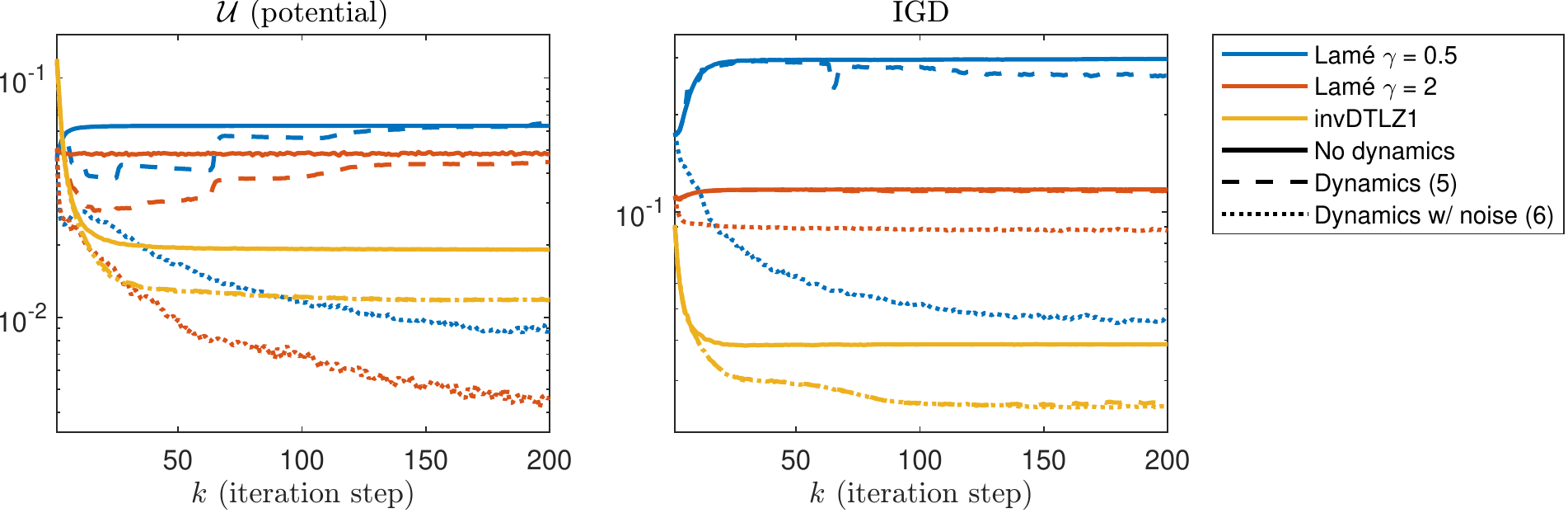}
\caption{Evolution of two performance metrics: the Morse potential and the IGD as functions of the iteratione step $k$ with different repulsive dynamics for the tri-objective problems considered.}
\label{fig:5}
\end{figure}

When testing the algorithm against tri-objective problems, we can only use dynamics \eqref{eq:new} and \eqref{eq:newnoise}. In all considered tests the proposed adaptive strategies improve the distribution of the computed solutions ($N=66$) over the front, even when the optimal distribution over the simplex differs significantly from the uniform one, as shown in Figs. \ref{fig:2}, \ref{fig:3} and \ref{fig:4}. When $m=3$ the improvement given by the stochastic component is evident in the corners of the front which are better approximated by the algorithm, see Figs. \ref{fig:3}c, \ref{fig:4}c. Quantitative estimates on the overall quality of the solutions are presented in Fig.\ref{fig:5} as functions of the iterative step $k$. The Morse potential measure $\mathcal{U}$ measures how well-distributed the solutions $\{F^i_k\}_{i=1}^N$ are over the front, while the Inverted Generational Distance (IGD) quantifies the distance between every points of a reference approximation and the set $\{F^i_k\}_{i=1}^N$ (see \cite{borghi2022adaptive} for a precise definition). For the Lam\'e problems considered, adding noise significantly improves the solution quality with respect to both measures, while for the inverted DTLZ1 problem the noise does not play a role. In all cases, the adaptive strategies allow to reach better solutions in terms of IGD. 

We remark that the computational cost per step $k$ is $\mathcal{O}(N^2)$, but it can simply be reduced by using random batch techniques \cite{AlPa,JLJ}, that is, evolving only a subset of parameters or estimating the forces with only a subset of point. We refer in particular to \cite{JLJ} for an example of random batch techniques applied to a system of charged particles over the sphere.


%
%
%
%
%
%
%
%
%
%
%
%

\section{Conclusions}

\label{s4}

We proposed an adaptive heuristic strategy to select scalarized sub-problems and obtain a uniform approximation of the Pareto front. Given that the sub-problems are parametrized by a vector of weights belonging to the unit simplex, the parameters are adapted following a repulsive dynamics where noise can also be added. We proposed different strategies that can be used for any number $m>1$ of objectives and make use of an arbitrary auxiliary method to solve the scalar sub-problems.
Numerical tests, both with $m=2$ and $m=3$ objectives, showed the validity of the repulsive dynamics to reach a better approximation of the Pareto front.

We observe that adding stochasticity improves the front approximation in some test problems considered. In future work, we plan on testing different stochastic dynamics which could be more effective than the vanilla Gaussian noise proposed in this work, exploiting, for instance, the the Aitchison geometry of the simplex \cite{lehmann2020darwinian}.

\section*{Acknowledgments} 
The work of G.B. is funded by the Deutsche Forschungsgemeinschaft (DFG, German Research Foundation) – Projektnummer 320021702/GRK2326 – Energy, Entropy, and Dissipative Dynamics (EDDy).

\end{document}